\documentclass[righttag,12pt]{amsart}
\usepackage{amssymb}
\usepackage{euscript}

\let\cal\mathcal

\newtheorem{prop}{Proposition}[section]
\newtheorem{theorem}[prop]{Theorem}
\newtheorem{cor}[prop]{Corollary}

\newtheorem{lem}[prop]{Lemma}

\theoremstyle{definition}
\newtheorem{defn}[prop]{Definition}
 
\theoremstyle{remark}
\newtheorem{rem}[prop]{Remark} 
\numberwithin{equation}{section}

\newcommand{\bbR}{{\mathbb{R}}}

\newcommand{\al}{\alpha}

\newcommand{\de}{\delta}
\newcommand{\e}{\varepsilon}

\newcommand{\la}{\lambda}

\newcommand{\ro}{\varrho}

\newcommand{\sgn}{\operatorname{sgn}}

\newcommand{\disp}{\displaystyle}

\newcommand{\lb}{\label}

\newcommand{\da}{\dot{a}}
\newcommand{\dda}{\ddot{a}}
\newcommand{\dx}{\dot{x}}
\newcommand{\ddx}{\ddot{x}}
\newcommand{\ddy}{\ddot{y}}

\newcommand{\no}{\nonumber}
\newcommand{\lra}{\longrightarrow}

\newcommand{\ONTO}{\buildrel {\mbox{\small onto}}\over \longrightarrow}

\begin{document}

\title{A pseudo-Daugavet property for narrow projections in Lorentz spaces}

\author{Mikhail M. Popov}  
\address{Department of Mathematics\\ 
Chervivtsi National University\\ str. Kotsjubyn'skogo 2, 
Chervivtsi, 58012 Ukraine} \email{popov@chv.ukrpack.net}

\author{Beata Randrianantoanina$^*$}\thanks{$^*$Participant, NSF Workshop
in Linear Analysis and Probability, Texas A\&M University}

\address{Department of Mathematics and Statistics
\\ Miami University \\ Oxford, OH 45056, USA}

 \email{randrib@muohio.edu}



\begin{abstract}
Let $X$ be a  rearrangement-invariant  space. An operator $T: 
X\to X$ is  called  narrow if for each measurable set $A$ and 
each $\varepsilon > 0$ there exists $x \in X$ with $x^2= \chi_A,\ 
\int x d \mu = 0$ and $\| Tx \| < \varepsilon$. In  particular all 
compact operators are narrow. We prove that if $X$ is a Lorentz 
function space $L_{w,p}$ on [0,1] with $p>2$, then there exists a 
constant $k_X>1$ so that for every narrow projection $P$ on $L_{w,p}$ 
$\| Id - P \| \geq k_X. $ This generalizes earlier results on 
$L_p$ and partially answers a question of E. M. Semenov. 
Moreover we prove  that every rearrangement-invariant  
function space $X$ with an 
absolutely continuous norm contains a complemented subspace 
isomorphic to $X$ which is the range of a narrow projection 
and a non-narrow projection, 
which gives a negative answer to a question of A.Plichko and M.Popov.
\end{abstract}
\maketitle

\section{Introduction}

We study a question of E. M. Semenov whether for every separable 
rearrangement-invariant function space $X \ (X \neq L_2)$ on 
$[0,1]$ there exists a constant $k_X>1$ so that for every rich 
subspace $Y \varsubsetneq X$ and every projection $P: X {\ONTO} 
Y$ we have
\begin{equation} \lb{S}
\| P \| \geq k_X.
\end{equation}

A subspace $Y \subset X$ is said to be rich if the quotient map 
$Q: X \to X/Y$ is narrow, where narrow operators are a 
generalization of compact operators on 
rearrangement-invariant function spaces (see Definitions~\ref{defn} 
and \ref{rich}). 
Narrow operators were formally introduced by Plichko and Popov 
\cite{PP90} but even before then they were used in the theory of 
Banach spaces, see e.g. \cite[Section 9]{JMST}. They were 
systematically studied in particular in \cite{PP90,KP96,KSW}.

 Natural examples of rich subspaces include 
subspaces of finite co-dimension, see \cite[\S 10]{PP90} for a 
discussion of properties of rich subspaces, here we just mention 
that if $Y \subset X$ is rich and operator $T$ on $X$ is compact 
then the subspace $(I+T)(Y)\subset X$ is rich (here $I$ denotes 
the identity operator on $X$) \cite[Proposition~10.11]{PP90}.

Results related to the above question of E. M. Semenov go back to 
 Lozanovski{\u\i} \cite{L66} and Benyamini and 
Lin \cite{BL85} who proved that for every $p, \ 1\leq p < \infty,\ p 
\neq 2$, there exists a function $\varphi_p:(0,\infty)\to 
(0,\infty)$ so that for every nonzero compact operator $T$ on 
$L_p[0,1]$ we have
\begin{equation}\lb{PD}
\| I - T \| \geq 1 + \varphi_p (\| T \|).
\end{equation}

When $p=1,\ \infty$, we can take $\varphi_p(t)=t$ and the above 
statement is called the {\it Daugavet equation}, cf. e.g. 
\cite{KSSW00}.  If $p=2$ we may have $\|I - T \| = 1$ for many 
compact operators $T$, e.g. for all orthogonal projections of 
finite rank.  If $1<p<\infty,\ p \neq 2$, then $\varphi_p(t)< t$ 
for all $t \in (0, \infty).$
Following \cite{O} we say that a Banach space $X$ satisfies a 
{\it pseudo-Daugavet property} if there exists a function 
$\varphi_X:(0,\infty)\to (0,\infty)$ so that \eqref{PD} is satisfied
for every nonzero 
compact operator $T$ on $X$.
Both the Daugavet property and the pseudo-Daugavet property
have important geometric implications, for example spaces
with the Daugavet property cannot be reflexive and  
the pseudo-Daugavet property is related
to the problem of best compact approximation in $X$,
i.e. to the question
whether every
bounded linear operator on $X$, has an element of best
approximation in the class of compact operators on $X$, see \cite{ABJS79}
for the thorough  introduction of this subject and references
(cf. also \cite{BL85,BO87}). 

 Popov \cite{Popov87} proved that \eqref{S} is valid in 
$L_p,\ 1 \leq p<\infty,\ p \neq 2$, i.e. he proved that for each 
$p,\ 1 \leq p< \infty,\ p \neq 2$, there exists a constant $k_p > 
1$ so that for every rich subspace $Y \varsubsetneq L_p$ and 
every projection $P:L_p {\ONTO}Y$ we have 
\begin{equation}\lb{PP}
\| P \| \geq k_p.
\end{equation}
This result has  very important applications to the study of the
geometric structure of spaces $L_p$, in particular it follows from
\eqref{PP} that every ``well'' complemented (i.e. with constant of 
complementation smaller than $\max\{k_p,k_{p'}\}$, where $1/p+1/p'=1$) 
subspace of $L_p$ is isomorphic
to $L_p$ \cite{PP90}.

Franchetti \cite{F92} found the exact value of the constant 
$k_p$. He proved that $k_p = \| I-A \|_p$, where $A$ is the rank 
one projection defined by: 
\begin{equation}\lb{defA}
Ax \overset{ \text{def} }{=}\left(\int_\Omega x(s)d 
\mu(s)\right)\cdot \bold1
\end{equation}
(Note here that the projection $A$ is well defined in any r.i. 
space with finite measure.)
The norm of $\|I-A\|_p$ for $p\in(1,\infty)$ has been evaluated 
by Franchetti \cite{F90} and, independently, by Oskolkov (unpublished).

Plichko and Popov \cite[Theorem~9.7]{PP90} generalized \eqref{PD} 
for all narrow operators on $L_p,\ 1\leq p < \infty,\ p \neq 2$.  
Later Oikhberg \cite{O} proved \eqref{PD} for compact operators 
on non-commutative $L_p, \ 1<p<\infty,\ p \neq 2.$  The first 
results related to \eqref{S} and \eqref{PD} (but not the Daugavet 
equation) for spaces other than $L_p$ seem to be the following:

\begin{theorem} \cite[Theorem~4.3]{KR}
Suppose that $X$ is a separable real order-continuous K\"othe 
function space on $(\Omega, \mu),$ where $\mu$ is nonatomic and 
finite. 
 Let $Y \subset X$ be a subspace of codimension one and $P$ be any  
 projection from $X$ onto $Y$.  Then 
$$
\| P \| >1,
$$
unless $X$ contains a band isometric to $L_2$, i.e. unless there 
exists a set $B \subset \Omega$ and a nonnegative measurable 
function $w$ with $supp \ w = B$ so that for any $x \in X$ with 
$supp \ x \subset B$
$$
\| x \|_X = (\int |x|^2 w d \mu)^{\frac{1}{2}}.
$$
\end{theorem}

\begin{theorem}\cite[Theorem~1]{FS95}\lb{FS}
Let $X$ be a real r.i. space on $(\Omega, \Sigma, \mu)$ where 
$\mu$ is nonatomic and $\mu(\Omega)=1$.  If $X$ is not 
isometrically isomorphic to $L_2$ then
$$
\| I-A \|_X > 1,
$$
where $A$ is the rank one projection defined in \eqref{defA}.
\end{theorem}

\begin{theorem}\cite[Theorem~4]{pams}\lb{pams}
Let $X$ be a separable real nonatomic r.i. space on $[0,1]$ which 
is not isometric to $L_2$. Let $Y \varsubsetneq X$ be any 
subspace of finite codimension and $P$ be any projection from $X$ 
onto $Y$.  Then
$$
\| P \| > 1.
$$
\end{theorem}

E. M. Semenov also communicated to us that he has proved 
\eqref{S} for every r.i. space $X$ on $[0,1]\ (X \neq L_2)$ and 
every projection $P$ from $X$ onto a rich subspace $Y 
\varsubsetneq X$ with the additional condition that 
$P(\chi_{[0,1]})=0$ (unpublished).

In the present paper we prove that \eqref{S} is valid if $X$ is a 
Lorentz space $L_{p,w}[0,1]$ with $p>2$ (there are no 
restrictions on $w$), see Theorem~\ref{main}.

Our result is valid in both complex and real case.  We have not 
attempted to find the exact value of $k_X$.  E. M. Semenov 
suggested that the result of Franchetti mentioned above (that 
$k_p = \|I - A \|_p$) should generalize to all r.i. spaces $X$.

In the final section we study duals of narrow operators.
In \cite{PP90} it was proved that in general the 
conjugate operator $T^*$ to a narrow operator $T:E\lra E$ 
need not be narrow (for any r.i. space $E$ with $E^*$ having an 
absolutely continuous norm to consider the notion of narrow 
operators). This naturally leads us to study properties
of operators which are conjugate to narrow operators. 
We call such operators {\it *-narrow} (see Definition~\ref{sn}).
We prove that *-narrowness of an 
operator on a reflexive r.i. space 
is a property of the image under it of the unit ball (Proposition~\ref{2}). 
However we show that
 the notion of narrow operators cannot be formulated in terms 
of the image. Namely we prove that if
 $E$ is an r.i. space with an absolutely continuous 
norm, then there exists a complemented subspace $E_0$ of $E$ 
isomorphic to $E$ and for which there are two projections onto 
$E_0$, one of which is narrow and the second is not narrow
(Theorem~\ref{4}). This answers (negatively) a 
question posed by A. Plichko and M. Popov \cite[Question 2, p. 71]{PP90}.

\section{Preliminaries}

Let us suppose that $\Omega$ is a Polish space and that $\mu$ is 
a $\sigma-$finite Borel measure on $\Omega.$ We use the term 
K\"othe space in the sense of \cite[p. 28]{LT2}.  Thus a {\it 
K\"othe function space} $X$ on $(\Omega,\mu)$ is a Banach space 
of (equivalence classes of) locally integrable Borel functions 
$f$ on $\Omega$ such that:  \newline (1) If $|f|\le |g|$ a.e. and 
$g\in X$ then $f\in X$ with $\|f\|_X \le \|g\|_X.$\newline (2) If 
$A$ is a Borel set of finite measure then $\chi_A\in X.$

We say that $X$ is {\it order-continuous} if whenever $f_n\in X$ 
with $f_n\downarrow 0$ a.e. then $\|f_n\|_X\downarrow 0.$ $X$ has 
{\it the Fatou property} if whenever $0\le f_n\in X$ with 
$\sup\|f_n\|_X<\infty$ and $f_n\uparrow f$ a.e. then $f\in X$ 
with $\|f\|_X=\sup\|f_n\|_X.$

{\it A rearrangement-invariant function space (r.i. space)} is a 
K\"othe function space on $([0,1],\mu)$ where $\mu$ is the 
Lebesgue measure which satisfies the conditions:\newline (1) 
Either $X$ is order-continuous or $X$ has the Fatou property.  
\newline (2) If $\tau:[0,1]\to [0,1]$ is any measure-preserving 
invertible Borel automorphism then $f\in X$ if and only if 
$f\circ\tau\in X$ and $\|f\|_X=\|f\circ\tau\|_X.$ \newline (3) 
$\|\chi_{[0,1]}\|_X=1.$

Next we recall the definition of the Lorentz spaces. These were 
introduced by Lorentz \cite{Lor50,Lor51}  in connection with some 
problems of harmonic analysis and interpolation theory. Since 
then they were extensively studied by many authors.
 
 If $f$\ is a measurable 
function, we define the {\it non-increasing rearrangement} of 
$f$\ to be
$$ f\sp *(t) = \inf\{\, s : \mu(| f| > s) \le t \,\}
.$$ 

Notice that when $f$ is a simple function, $ f = \sum^m_{k=1} a_k 
\chi_{A_k}$, then $f^*$ is also a simple function and the range of
$f^*$ equals $\{|a_k|\ : \ k=1,\dots,m\}$.

If $1\le p < \infty$, and if $w:(0,1)\to(0,\infty)$\ is a 
non-increasing function, we define the {\it Lorentz norm} of a 
measurable function $f$ to be
$$ \| f\|_{w,p} = \left(\int_{[0,1]} w(t) f\sp *(t)\sp p
\,dt\right)\sp {1/p} .$$ We define the {\it Lorentz space} $L\sb 
{w,p}([0,1],\mu)$ to be the space of those measurable functions 
$f$\ for which $\|f\|\sb {w,p}$\ is finite. These spaces are a 
generalization of the $L\sb p$\ spaces: if $w(x) = 1$\ for all 
$0\le x <1$, then $L\sb {w,p} = L\sb p$\ with equality of norms. 
Lorentz spaces are one of the most important examples of r.i. 
spaces.

Narrow operators generalize the notion of compact operators on 
r.i. spaces.  
\begin{defn}\lb{defn}
Let $X$ be an r.i. space and $Y$ be any Banach space.  We say 
that an
 operator $T: X\to Y$ is {\it narrow} if for each measurable set $A$ 
 and each $\e > 0$ there exists $x \in X$ with 
 $x^2= \chi_A,\ \int x d \mu = 0$ and $\| Tx \| < \e$.
\end{defn}

In fact in the above definition the condition that $\int x d \mu 
=0$ can be omitted \cite[Proposition~8.1]{PP90}. Every compact 
operator is narrow \cite[Proposition~8.2]{PP90}.

Plichko and Popov proved the following lemma which will be useful 
in our proof.

\begin{lem}\lb{pp}\cite[Lemma~8.1]{PP90}
Let $T: X \to Y$ be narrow.  Then for every $\e>0$, every 
measurable set $A$ and every integer $n \geq 1$ there exists a 
partition $A= A' \cup A''$ into measurable subsets with $\mu(A')= 
2^{-n}\mu (A)$ and $\mu(A'')=(1-2^{-n}) \mu(A)$ such that $\|Th 
\| < \e,$ where $h =(2^n -1)\chi_{A'}-\chi_{A''}$.
\end{lem}

\begin{defn}\lb{rich}
Let $X$ be an r.i. space.  A subspace $Y \subset X$ is called 
{\it rich} if the quotient map $Q: X \to X/Y$ is narrow.

In other words, $Y \subset X$ is rich if for every measurable set 
$A$ and every $\e > 0$ there exist $y \in Y$ and $x \in X$ so 
that $x^2= \chi_A,\ \int x d \mu =0$ and $\| x-y \| < \e$.
\end{defn}

In particular every subspace of finite codimension is rich.

\section{Main result}

\begin{theorem}\lb{main}
Suppose $L_{p,w}$ is a Lorentz space on $[0,1]$ with $p>2$.  Then 
there exists $\ro_p >1$ so that for every nontrivial projection 
$P$ from $L_{p,w}$ onto a rich subspace
$$
\| P \| \geq \ro_p.
$$
\end{theorem}

In the proof of Theorem~\ref{main} we will use the following two 
propositions:

\begin{prop} \lb{case1}
Suppose $L_{p,w}$ is a Lorentz space on $[0,1]$ with $p>2$.  Then 
there exist $\de_p \in (0,{1}/{8}), \la_p=\la_p(\de_p, p) \in 
({\de_p}/({\de_p -4}), 0)$ and $\gamma_p = \gamma_p(\lambda_p, 
\de_p, p) \in (0,1)$ so that
$$
\gamma_p + \frac{1}{2}| \lambda_p | \de_p < 1
$$
which satisfy the following property: 

For every simple function $x = \sum^m_{k=1}a_k \chi_{A_k}$ so 
that $1 \leq \| x \|_{p,w} \leq 1 + \frac{3}{2}\de_p$ and
$$
\frac{|a_i|}{|a_j|} \notin (3-\de_p, 3)
$$
for all $i, j=1, \ldots, m$; and for every partition $A_k = B_k 
\sqcup C_k$ with $\mu(B_k)=({1}/{4})\mu(A_k)$ we have
$$
\big\| \lambda_px + \sum^m_{k=1}a_k(3\chi_{B_k}-\chi_{C_k}) 
\big\|_{p,w} \leq \gamma_p \big\| \sum^m_{k=1} a_k 
(3\chi_{B_k}-\chi_{C_k}) \big\|_{p,w}.
$$
\end{prop}

\begin{prop}\lb{ddot}
Let $X$ be an r.i. space. Given a simple function $x = 
\sum^m_{k=1}a_k \chi_{A_k}\in X$ and $\de \in (0,1/8)$ there 
exists a simple function $\ddx=\ddx(\de)=\sum^m_{k=1}\dda_k 
\chi_{A_k}$ so that for all $i, j=1, \ldots,m$,
\begin{equation*}
\frac{|\dda_i|}{|\dda_j|} \notin (3-\de,3)
\end{equation*}
and $\| x - \ddx \|< (3/2)\de$, 
$\|x\|\le\|\ddx\|<(1+(3/2)\de)\|x\|$.
\end{prop}

Let us first show that Theorem~\ref{main} is indeed a consequence 
of Propositions~\ref{case1} and \ref{ddot}.

\begin{proof}[Proof of Theorem~\ref{main}]

Fix $\e > 0$.  Since $P$ is a non-trivial projection, there 
exists a simple function $x=\sum^m_{k=1}a_k \chi_{A_k},\ (a_k 
\neq 0)$ with $\| x \| =1$ and $\|P x \| < \e$ (note that since 
we will always work in $L_{p,w}$ we will drop the subscript and 
simply use $\| \cdot \|$ to mean $\| \cdot \|_{p,w}$ throughout 
this proof).

Let $\de_p$ be as defined in the statement of 
Proposition~\ref{case1}.  Since $\de_p \in (0,{1}/{8})$, by 
Proposition~\ref{ddot}, there exists a simple function $\ddx = 
\sum^m_{k=1}\dda_k \chi_{A_k}$ with $1 \leq \| \ddx \| < 1 + 
({3}/{2}) \de_p,\ \| x-\ddx \| < (3/2)\de_p$ and so that
\begin{equation}\lb{dd3}
\frac{| \dda_i |}{| \dda_j |} \notin (3-\de_p, \de_p),
\end{equation}
for all $i, j=1, \ldots, m$.

Since $I-P$ is narrow, by Lemma~\ref{pp}, for each $k,\ 1 \leq k 
\leq m$, there exists a partition $A_k = B_k \sqcup C_k$ so that 
$\mu(B_k)= ({1}/{4})\mu(A_k)$ and 
$$
\|(I-P)(3\chi_{B_k}-\chi_{C_k}) \| < \frac {\e}{| \dda_k |m}.
$$
Then for $\ddy = \sum^m_{k=1} \dda_k (3 \chi_{B_k}-\chi_{C_k})$  
we obtain
$$
\|(I-P)\ddy \| < \e,
$$
and
$$
\| \ddy \| \leq 3 \| \ddx \| < 3 + \frac{9}{2} \de_p.
$$

Moreover, by \eqref{dd3} and Proposition~\ref{case1} we conclude 
that
$$
\| \lambda_p \ddx + \ddy \| \leq \gamma_p \| \ddy \|,
$$
where $\lambda_p$ and $\gamma_p$ are constants defined in 
Proposition~\ref{case1}. Thus
\begin{equation*}
\begin{split}
\| \ddy \| &=\| P(\lambda_p \ddx + \ddy)- P \lambda_p \ddx
- P\lambda_p x + P\lambda_p x + \ddy - P \ddy \| \\
&\leq \| P \| \cdot \| \lambda_p \ddx + \ddy \| + | \lambda_p | 
\cdot \| P \| \cdot \|
\ddx -x \| + | \lambda_p | \|Px \| + \| (I-P) \ddy \| \\
&\leq \| P \| \cdot \gamma_p \| \ddy \| + \| P \| |\lambda_p| 
\cdot \frac{3}{2}\de_p + 
| \lambda_p | \e + \e \\
&\leq \| \ddy \| \cdot \| P \| \left(\gamma_p + \frac{3\de_p 
|\lambda_p|}{6 + 9 \de_p}\right) + \e
(| \lambda_p |+1)\\
&\leq \| \ddy \| \cdot \| P \|\left(\gamma_p + \frac{1}{2} \de_p 
| \lambda_p|\right) + \e (|\lambda_p|+1).
\end{split}
\end{equation*}

Since $\e$ was arbitrary we obtain 
$$
\| P \| \geq (\gamma_p+ \frac{1}{2}\de_p | \lambda_p |)^{-1} 
\overset{ \text{def} }{=} \varrho_p.
$$

By Proposition~\ref{case1}, $\varrho_p >1$.
\end{proof}

\begin{rem} Note that the same proof will demonstrate that 
whenever $T$ is a narrow operator on $L_{p,w}$, $p>2$, such that 
1 is an eigenvalue of $T$, i.e. such that there exists a nonzero 
element $x\in L_{p,w}$ with $Tx=x$, then 
$$\|I-T\|\ge \varrho_p 
>1.$$
(Simply replace $P$ in the proof by $I-T$, and note that 
Propositions~\ref{case1} and \ref{ddot} do not depend on the 
operator at all.)
\end{rem}

\begin{proof}[Proof of Proposition~\ref{case1}]
Let $\de \in (0, {1}/{8})$ and $x = \sum^m_{k=1} a_k \chi_{A_k}$ 
be a simple function so that $1 \leq \| x \|_{p,w} \leq 1 + 
({3}/{2})\de$ and
\begin{equation}\lb{not3}
\frac{|a_i|}{|a_j|} \notin (3-\de,3),
\end{equation}
for all $i, j=1, \ldots, m$.  We assume without loss of 
generality that $|a_1| \geq |a_2| \geq \ldots \geq |a_m|$.  Let 
$B_k, C_k$ for $k=1, \ldots, m$, be subsets as described in the 
statement of the proposition.

Denote
$$
y = \sum^m_{k=1} 3a_k \chi_{B_k}-a_k \chi_{C_k},
$$
and set $b_k = 3a_k, c_k = -a_k$ for $k=1,\ldots,m$.

Thus
\begin{equation*}
\begin{split}
y &= \sum^m_{k=1}b_k \chi_{B_k} + c_k \chi_{C_k}, \\
\lambda x + y &=\sum^m_{k=1} b_k \big(1+\frac{\lambda}{3}\big) 
\chi_{B_k} + c_k (1-\lambda) \chi_{C_k}.
\end{split}
\end{equation*}

We first notice that if $0 \geq \lambda > {\de}/({\de-4})>-1$ 
then for all $i,j=1, \ldots, m$ the following hold:
\begin{eqnarray}
|b_i|\big(1+\frac{\lambda}{3}\big) &\leq& |b_j|\big(1+\frac{\lambda}{3}\big)\Longleftrightarrow |b_i| \leq |b_j|;  \lb{b}\\
|c_i| (1-\lambda) &\leq& |c_j|(1-\lambda) \Longleftrightarrow |c_i| \leq |c_j|;  \lb{c}  \\
|c_i| (1-\lambda) &<& |b_i|\big(1+\frac{\lambda}{3}\big) \ \ \ \text{for all \ $i$};   \lb{cb}  \\
|b_i| \big(1+\frac{\lambda}{3}\big) &\leq& |c_j|(1-\lambda) 
\Longleftrightarrow
|b_i| \leq |c_j|;  \lb{bc} \\
|b_i| \big(1+ \frac{\lambda}{3}\big) &\neq& |c_j| (1-\lambda)  \ 
\ \ \text{for all \ $i,j$}. \lb{noneq}
\end{eqnarray}

Indeed \eqref{b}, \eqref{c} and \eqref{cb} are obvious since 
$1+\lambda/3 >0$, $1-\lambda >0$ and $1-\lambda < 3(1+ 
\lambda/3)$. To see \eqref{bc}``$\Rightarrow$'', suppose for 
contradiction that there exist $i,j$ so that
$$
|b_i|(1+\frac{\lambda}{3}) \leq |c_j| (1-\lambda),\ \ \text{and}\ 
\ |b_i|>|c_j|.
$$
Then
$$
1>\frac{|c_j|}{|b_i|}\geq \frac{1-\lambda}{1+\frac{\lambda}{3}} = 
\frac{1-\lambda}{3(3+\lambda)}.
$$
Thus, since $\lambda \in({\de}/({\de-4}),0)$,
$$
3>\frac{|a_j|}{|a_i|}\geq \frac{1-\lambda}{3+\lambda}
>\frac{3-\frac{\de}{4-\de}}{1+\frac{\de}{4-\de}} = 3-\de,
$$
which contradicts \eqref{not3} and \eqref{bc}``$\Rightarrow$'' is 
proved.

Next, suppose $|b_i|\leq|c_j|$. Since $\lambda<0$ we get 
$$
|b_i|(1+\frac{\lambda}{3})<|b_i|\leq|c_j|<|c_j|(1-\lambda).
$$

Thus \eqref{bc}``$\Leftarrow$'' and \eqref{noneq} are proved.

Now define numbers $t_{C_k},t_{B_k}$ for $k=1, \ldots,m$ as 
follows:
\begin{equation*}
\begin{split}
t_{C_i} &=\sum_{k<i} \mu(C_k) + \sum_{l:|b_l|>|c_i|}\mu(B_l) , \\
t_{B_j} &=\sum_{k<j} \mu(B_{k}) + \sum_{i:|c_i|\geq |b_j|} 
\mu(C_i).
\end{split}
\end{equation*}

It follows from \eqref{b}-\eqref{noneq} that for all 
$i,j=1,\dots,m$
\begin{equation}\lb{order}
\begin{split}
t_{C_i}>t_{C_j} &\ \ \Rightarrow \ \  \Big(|c_i|\leq|c_j|\ \ 
\text{and} \ \ |c_i|(1-\lambda)
\leq |c_j|(1-\lambda)\Big), \\
\Big(|c_i|<|c_j|, &\text{\ \ or, equivalently,\ \ } 
|c_i|(1-\lambda)<|c_j|(1-\lambda)\Big)
\ \ \Rightarrow \ \  t_{C_i}>t_{C_j}, \\
t_{B_i}>t_{B_j} &\ \ \Rightarrow \ \ \Big(|b_i|\leq|b_j| \text{\ 
\ and\ \ } |b_i|(1+\frac{\lambda}{3})
\leq |b_j|(1+\frac{\lambda}{3})\Big), \\
\Big(|b_i|<|b_j| &\text{\ \ or, equivalently,\ \ } 
|b_i|(1+\frac{\lambda}{3}) <|b_j|(1+\frac{\lambda}{3})\Big)
\ \ \Rightarrow \ \ t_{B_i}>t_{B_j}, \\
t_{B_i}>t_{C_i} &\ \ \Leftrightarrow \ \ 
|b_i|(1+\frac{\lambda}{3}) <|c_j|(1-\lambda)\ \  \Leftrightarrow 
\ \ |b_i|\leq |c_j|.
\end{split}
\end{equation}

Now define weights:
$$
w_{B_k}=\int^{t_{B_k}+\mu(B_k)}_{t_{B_k}} w d \mu \ , \ \ \ \ \ 
w_{C_k}=\int^{t_{C_k}+\mu(C_k)}_{t_{C_k}} w d \mu.
$$

By \eqref{order} we obtain
\begin{equation*}
\begin{split}
\|y\|^p_{p,w} &=\sum^m_{k=1}\big[|b_k|^p w_{B_k}+|c_k|^p w_{C_k}\big], \\
\| \lambda x+y \|^p_{p,w} &= 
\sum^m_{k=1}\big[|b_k|^p(1+\frac{\lambda}{3})^p w_{B_k} + 
|c_k|^p(1-\lambda)^p w_{C_k}\big].
\end{split}
\end{equation*}

(Heuristically speaking, the non-increasing order of moduli of 
coefficients of $y$ is the same as the  non-increasing order of 
moduli of coefficients of $\lambda x + y$).

Thus if we set
$$
\psi(\lambda)=\sum^m_{k=1}\big[|b_k|^p(1+\frac{\lambda}{3})^p 
w_{B_k} +|c_k|^p(1-\lambda)^p w_{C_k}\big],
$$
for $\lambda \in(-3,1)$, then
\begin{equation}
\begin{split}
\psi(0) &= \| y\|^p_{p,w}, \\
\psi(\lambda) &= \| \lambda x + y \|^p_{p,w} \ \ \text{for} \ \  
\lambda \in (\de/(\de-4), 0).
\end{split}
\end{equation}

Clearly $\psi$ is differentiable for all $\lambda \in (-3,1)$ and
$$
\psi'(\lambda)=\sum^m_{k=1}\big[p|b_k|^p(1+\frac{\lambda}{3})^{p-1} 
\frac{1}{3} w_{B_k}- p |c_k|^p (1-\lambda)^{p-1} w_{C_k}\big].
$$

Thus
\begin{equation}\lb{devpsi}
\begin{split}
\psi'(0) &=p\sum^m_{k=1}\big[\frac{1}{3}|b_k|^p w_{B_k}-|c_k|^p w_{C_k}\big]\\
&=p \sum^m_{k=1} |a_k|^p \big[3^{p-1} w_{B_k}- w_{C_k}\big]
\end{split}
\end{equation}

We now need to compare the quantities $w_{B_k}$ and $w_{C_k}$ for 
a given $k,\ 1 \leq k \leq m$.  It follows from \eqref{cb} that 
$t_{C_k}>t_{B_k}$. Moreover, by definition of $B_k$ and $C_k$ we 
have $\mu(C_k)=3\mu(B_k)$.  Thus, since $w$ is non-increasing, we 
get
\begin{equation*}
\begin{split}
w_{C_k} &=\int_{t_{C_k}}^{t_{C_k}+\mu(C_k)} w d \mu =
\int_{t_{C_k}}^{t_{C_k}+ 3\mu(B_k)} w d \mu \\
&\leq \int_{t_{B_k}}^{t_{B_k}+\mu(B_k)} w d \mu + 
\int_{t_{B_k}+\mu(B_k)}^{t_{B_k}+2\mu(B_k)} w d \mu +
\int_{t_{B_k}+ 2\mu(B_k)}^{t_{B_k}+3\mu(B_k)} w d \mu \\
&\leq 3\int_{t_{B_k}}^{t_{B_k}+\mu(B_k)} w d \mu\\
&=3 w_{B_k}.
\end{split}
\end{equation*}

Thus for all $k=1, \dots,m,$
\begin{equation}\lb{est1}
3^{p-1}w_{B_k}-w_{C_k} \geq 3^{p-1} {w_{B_k}}-3 w_{B_k} = 
w_{B_k}(3^{p-1}-3).
\end{equation}

In analogy to numbers $t_{C_k},t_{B_k},w_{C_k},w_{B_k}$ we define:
$$
t_{A_k} =\sum_{l<k} \mu(A_l), 
$$
$$
w_{A_k} =\int_{t_{A_k}}^{t_{A_k}+ \mu(A_k)} w d\mu.
$$

Since we assumed that $|a_1| \geq |a_2|\geq \ldots \geq |a_m|$ we 
obtain
$$
\|x\|^p_{p,w}=\sum^m_{k=1} |a_k|^p w_{A_k}.
$$

Further, since for all $k, \ \mu(A_k) = \mu(B_k)+ \mu(C_k)$ and 
since $|c_l| \geq |b_j| \Rightarrow |c_l| >|c_j| \Rightarrow 
l<j$, we obtain
\begin{eqnarray*}
t_{B_j} &=& \sum_{k<j} \mu(B_k)+ \sum_{l:|c_l|\geq|b_j|} 
\mu{(C_l)} 
\leq \sum_{k<j} \mu(B_k)+\sum_{l<j} \mu(C_l)\\
&=& \sum_{k<j} \mu(A_k)= t_{A_j}.
\end{eqnarray*}

Therefore, since for all $j=1, \ldots, m,\ \mu(B_j)=\frac{1}{4} 
\mu(A_j)$, we obtain:
\begin{eqnarray*}
w_{B_j} &=& \int_{t_{B_j}}^{{t_{B_j}}+ \mu(B_j)}
 w d \mu
\geq \int_{t_{A_j}}^{{t_{A_j}}+\frac{1}{4}\mu(A_j)} w d \mu \\
&\geq& \frac{1}{4} \int_{t_{A_j}}^{{t_{A_j}} + \mu(A_j)}w d \mu = 
\frac{1}{4}w_{A_j}.
\end{eqnarray*}

Thus we can continue the estimate from \eqref{est1} as follows
$$
3^{p-1}w_{B_k}-w_{C_k} \geq w_{B_k}(3^{p-1}-3) \geq 
\frac{1}{4}(3^{p-1}-3) w_{A_k}.
$$

Plugging this into \eqref{devpsi} we get
\begin{equation}\lb{Cp}
\begin{split}
\psi'(0) &= p \sum^m_{k=1} |a_k|^p \big[3^{p-1} w_{B_k}-w_{C_k}\big]\\
&\geq \frac{1}{4}p (3^{p-1}-3) \sum^m_{k=1} |a_k|^p w_{A_k} \\
&= \frac{1}{4}p(3^{p-1}-3) \|x \|^p_{w,p}\\
&\geq \frac{1}{4}p(3^{p-1}-3)\overset{ \text{def }}{=} C_p.
\end{split}
\end{equation}

Note that our assumption that $p>2$ guarantees that $C_p>0$.

 Our next step is to estimate from above the value of 
$|\psi''(\lambda)|$ when $\lambda \in (\de/(\de-4), 
\de/(4-\de))$. We have, since $\de \in(0,\frac{1}{8})$,
\begin{equation}\lb{Mp}
\begin{split}
|\psi''(\lambda)| &=\Big|p(p-1)\sum^m_{k=1} 
\big[\frac{1}{9}|b_k|^p(1+\frac{\lambda}{3})^{p-2} 
w_{B_k} + |c_k|^p(1-\lambda)^{p-2}w_{C_k}\big]\Big| \\
&\leq p(p-1)(1+\frac{\de}{4-\de})^{p-2} \sum^m_{k=1}
\big[\frac{1}{9} |b_k|^p w_{B_k}+ |c_k|^p w_{C_k}\big]\\
&\leq p(p-1)(1+\frac{\de}{4-\de})^{p-2} \|y\|^p_{w,p}\\
&\leq p(p-1)(\frac{4}{4-\de})^{p-2}3^p \|x \|^p_{w,p}\\
&\leq p(p-1)(\frac{4}{4-\de})^{p-2} 3^p (1+\frac{3}{2}\de)^p\\
&\leq p(p-1)4^p  \overset{\text{def} }{=} M_p.
\end{split}
\end{equation}

By Taylor's Theorem for $\lambda \in ({\de}/({\de-4}),0)$ we get 
$$
\psi(\lambda)= \psi(0) + \lambda \psi'(0) + \frac{1}{2}\lambda^2 
\psi'' (\theta),
$$
where $\theta \in(\lambda,0)$. Note that: 
$$
\psi(0)= \| y \|^p_{p,w} \leq 3^p \|x|^p_{p,w} \leq 
[3(1+\frac{3}{2}\de)]^p =(3 + \frac{9}{2}\de)^p.
$$
 Thus by 
\eqref{Cp} and \eqref{Mp} we have 
\begin{eqnarray*}
\psi(\lambda) &\leq& \psi(0)+ \lambda C_p + \frac{1}{2} \lambda^2 M_p \\
&\leq& \psi(0)[1+\frac{\lambda}{(3+\frac{9}{2}\de)^p}(C_p + \frac{1}{2}\lambda M_p)] \\
&\leq& \psi(0)[1 +\lambda 3^{-p}(C_p +\frac{1}{2} \lambda M_p)].
\end{eqnarray*}

Thus when $|\lambda| \leq {C_p}/{M_p}$ we get
$$
\psi(\lambda) \leq \psi(0)\big[1 + \lambda \frac{C_p}{2 \cdot 
3^p}\big].
$$

If $\de \leq \min \{{1}/{8}, {4C_p}/{M_p}\}$ we set $\lambda = 
-{\de}/{4}$ and
$$
\gamma(\de)\overset{ \text{def} }{=}(1-\de \cdot \frac{C_p}{8 
\cdot 3^p})^{\frac{1}{p}}.
$$
Then $\gamma(\de)<1$ and by definition of the function $\psi$ we 
have
$$
\big\| - \frac{\de}{4}\cdot x+y \big\| \leq \gamma (\de)\| y \|.
$$

To finish the proof of the proposition we notice that by the 
Bernoulli inequality
$$
\gamma(\de)< 1 - \de \cdot \frac{C_p}{8 \cdot 3^p \cdot p}\ .
$$

Set
$$
D_p \overset{ \text{def} }{=}\frac{C_p}{8 \cdot 3^p \cdot p} = 
\frac{3^{p-2}-1}{32 \cdot 3^{p-1}}\ .
$$

Since $p>2$ we have $D_p>0$ and we obtain:
$$
\gamma(\de) + \frac{1}{2} |\lambda| \de = \gamma(\de)+ 
\frac{1}{8}\de^2 < 1-\de D_p + \frac{1}{8}\de^2 =1-\de(D_p - 
\frac{1}{8}\de).
$$
Thus if $\de \leq {8}D_p$ then
$$
\gamma(\de)+ \frac{1}{2} |\lambda| \de < 1.
$$

Hence we can take $\de_p = \min\{{1}/{8},{4C_p}/{M_p},8{D_p}\},\ 
\lambda_p = -{\de_p}/{4},\ \gamma_p = \gamma(\de_p) = 
(1-\de_ppD_p)^{\frac{1}{p}}$ and the proposition is proved.
\end{proof}

\begin{rem}\lb{p<2}
It is clear that the above proof does not work for $p<2$. Indeed, 
when $p<2$ then the estimate \eqref{est1} becomes meaningless and 
both constants $C_p$ and $D_p$ are negative. Moreover, for every 
$p$, $1\le p<2,$ it is not difficult to construct weights $w_p$ 
so that when $x=\chi_{[0,1]}$ is partitioned into any disjoint 
sets
 $[0,1] = B \sqcup C$ with $\mu(B)={1}/{4}$ then for any $\lambda\in\bbR$:
\begin{equation}\lb{eq}
\| \lambda x + (3\chi_{B}-\chi_{C}) \|_{p,w_p} \geq 
 \| 3\chi_{B}-\chi_{C} \|_{p,w_p}.
\end{equation}
In fact one can take e.g. 
$$ w_p= \frac4{3^{p-1}+1}\chi_{[0,\frac14)}+ 
\frac{4\cdot 3^{p-2}}{3^{p-1}+1}\chi_{[\frac14,1]}.$$
This is a well defined weight when $1\le p<2$. It is routine, even though 
tedious, to check that $L_{p, w_p}$ satisfy \eqref{eq} for all $p$ with
$1\le p<2$. We leave the details to the interested reader.

We suspect that if $1\le p<2$ the Proposition~\ref{case1} fails in $L_{p,w}$
for any nonconstant weight $w$, but we have not 
checked it carefully. This, of course, does not mean that we 
believe that Theorem~\ref{main} fails for $p<2$. For some comments
involving duality please see Section~\ref{dual}.
\end{rem}

\section{Proof of Proposition~\ref{ddot}}

The first step of the proof of Proposition~\ref{ddot} is the 
following lemma.

\begin{lem}\lb{step1}
Let $X$ be an r.i. space and $x \in X$ be a simple function $x = 
\sum^m_{k=1}a_k \chi_{A_{k}}$.  For any $\eta > 0$ there exists 
$\dx = \dx(\eta)=\sum^m_{k=1} \da_k \chi_{A_k}$ so that for all 
$i, j=1, \ldots, m$,
$$
\frac{|\da_{j}|}{|\da_{i}|} \notin (1, 1 + \eta)
$$
and $\|x - \dx \| < \eta \| x \|$, $\|\dx\|\ge\|x\|$.
\end{lem}

\begin{proof}[Proof of Lemma~\ref{step1}]
Without loss of generality we assume that $|a_1| \geq | a_2 | 
\geq \ldots \geq |a_m|$.

Let $r_0 = 1 < r_1 < r_2 \ldots < r_n = m$ be such that
\begin{equation*}
\frac{|a_j|}{|a_i|}
\begin{cases}
< 1 + \eta     &\text{if there exists $k$ with $r_k \leq j < i < r_{k+1}$,}\\
\geq 1 + \eta  &\text{if there exists $k$ with $j < r_k \leq i$.}
\end{cases}
\end{equation*}

Define 
$$
\da_j = \sgn (a_j)|a_{r_{k(j)}}|,
$$
where $k(j)$ is such that $r_{k(j)}\leq j < r_{k(j)+1}$.

Then for any $j < i$ we have $k(j)\le k(i)$ and
\begin{equation*}
\frac{|\da_j|}{|\da_i|} = \frac{|a_{r_{k(j)}}|}{|a_{r_{k(i)}}|}
\begin{cases}
= 1           &\text{if $k(j)= k(i)$,}\\
\geq 1+ \eta  &\text{if $k(j)< k(i)$.}
\end{cases}
\end{equation*}

Thus 
$$
\frac{|\da_{j}|}{|\da_{i}|} \notin (1, 1 + \eta)
$$
as required.

Moreover for all $j=1, \ldots, m:$
$$
\frac{\da_j}{a_j}= \frac{|a_{r_{k(j)}}|}{|a_j|} \in [1, 1 + \eta).
$$

Thus $\|\dx\|\ge\|x\|$ and
\begin{equation*}
\begin{split}
\|\dx - x\| &= \big\|\sum^m_{j=1}(\da_j - a_j) \chi_{A_{k}}\big\|
= \big\| \sum^m_{j=1} a_j (\frac{\da_j}{a_j}-1) \chi_{A_{k}}\big\|\\
&< \eta \big\| \sum^m_{j=1} a_j \chi_{A_j}\big\|=\eta \| x \|,
\end{split}
\end{equation*}
which ends the proof of the lemma.
\end{proof} 

In the next lemma we gather, for easy reference, a few simple 
arithmetic inequalities which will be useful in the proof of 
Proposition~\ref{ddot}.

\begin{lem}\lb{close}
Let $\de \in (0,3)$ and $t_i, t_i,$ $i=1,2,3,4,$ be positive real 
numbers such that
$$
\frac{\tilde{t_i}}{t_i} \in \left[1, \frac{3}{3-\de}\right) \ 
\text{for} \ i=1,2,3,4.
$$
Then we have:
\begin{itemize}
\item [$(i)$] 
If $t_i = t_2, \ \tilde{t_i}/\tilde{t_3} \in (3-\de,3), \ 
\tilde{t_2}/\tilde{t_4}  \in (3-\de,3)$, then  
$\tilde{t_3}/\tilde{t_4} \in ((3-\de)^{3}/27,{27}/{(3-\de)^{3}})$.
\item[$(ii)$]  
If ${t_1}/{t_2} \in (3-\de,3),\  {t_3}/{t_2}=3$, then 
${t_3}/{t_1} \in (1, {3}/{(3-\de)})$.
\item[$(iii)$]
If $t_1 < t_2$, then ${\tilde{t_1}}/{\tilde{t_2}} < 
{9}/{(3-\de)^2}$.
\item[$(iv)$]
If ${t_1}/{\tilde{t_2}}=3,\ {t_1}/{\tilde{t_3}} \in (3-\de,3)$, 
then ${t_2}/{t_3} \in ({(3-\de)^2}/{9}, {3}/{(3-\de)})$.
\end{itemize}
\end{lem}

\begin{proof}[Proof of Lemma~\ref{close}]
The proofs of Lemma~\ref{close}$(i)-(iv)$ are very simple and 
very similar to each other.  As an illustration we prove the 
implication $(i)$.

We have
\begin{equation*}
\frac{t_3}{t_4}=\frac{t_3}{\tilde{t_3}} \cdot \frac{\tilde{t_3}} 
{\tilde{t_1}}\cdot \frac{\tilde{t_1}}{t_1} \cdot \frac{t_1}{t_2} 
\cdot \frac{t_2}{\tilde{t_2}} \cdot 
\frac{\tilde{t_2}}{\tilde{t_4}} \cdot \frac{\tilde{t_4}}{t_4}
\end{equation*}

Since
\begin{equation*}
\begin{split}
\frac{t_3}{\tilde{t}_3} &\in \left(\frac{3-\de}{3},1\right],\ 
\frac{\tilde{t_3}} {\tilde{t}_1} \in 
\left(\frac{1}{3},\frac{1}{3-\de}\right),\ 
\frac{\tilde{t_1}}{t_1} \in \left[1, 
\frac{3}{3-\de}\right), \ \frac{t_1}{t_2}=1,\\
 \frac{t_2}{\tilde{t_2}} 
&\in \left(\frac{3-\de}{3},1\right],\ \frac{\tilde{t_2}} 
{\tilde{t_4}} \in(3-\de,3),\ \frac{\tilde{t_4}}{t_4} \in 
\left[1,\frac{3}{3-\de}\right),
\end{split}
\end{equation*}
we obtain
\begin{equation*}
\begin{split}
\frac{t_3}{t_4} &\in \left(\frac{3-\de}{3}\cdot \frac{1}{3}\cdot 
1 \cdot 1 \cdot \frac{3-\de}{3}\cdot (3-\de) \cdot 1, 1 \cdot 
\frac{1}{3-\de}\cdot \frac{3}{3-\de} \cdot 1 \cdot 1 \cdot 3 
\cdot \frac{3}{3-\de}\right) \\
&= \left(\frac{(3-\de)^3}{27},\frac{27}{(3-\de)^3}\right). 
\end{split}
\end{equation*}

Implications $(ii)$-$(iv)$ are proved in a very similar way.
\end{proof}

In the proof of Proposition~\ref{ddot} we will also need the 
following definition:
\begin{defn}
For any simple function $y = \sum^m_{k=1}d_k \chi_{D_k}$  define 
sets:
\begin{equation*}
S(y,k)=\{j \in \{1,\ldots,m\}: \frac{|d_j|}{|d_k|}\in (3- 
\de,3)\}.
\end{equation*}
\end{defn}

\begin{proof}[Proof of Proposition~\ref{ddot}]
Let $\eta>0$ be such that 
$$
1+ \eta = \left(\frac{3}{3-\de}\right)^3.
$$

By Lemma~\ref{step1} there exists 
$\dx=\dx(\eta)=\sum^m_{k=1}\da_k \chi_{A_k}$ with
\begin{equation}\lb{xdot}
\| x-\dx \| < \eta \| x \|,
\end{equation}
so that for all $i,j = 1, \ldots, m$,
\begin{equation*}
\frac{| \da_i |}{| \da_j |} \notin (1,1+\eta).
\end{equation*}

By symmetry this means that for all $i,j=1,\ldots,m$,
\begin{equation}\lb{eta}
\frac{|\da_i|}{|\da_j|} \in \left(\frac{1}{1+\eta},1+\eta\right) 
\Longrightarrow |\da_i | = |\da_j|.
\end{equation}

To prove the lemma we need to construct a simple function $\ddx$ 
so that $\| x - \ddx \|< (3/2)\de$, 
$\|x\|\le\|\ddx\|<(1+(3/2)\de)\|x\|$ and
\begin{equation}\lb{xddot}
S(\ddx,k)= \emptyset \ \text{for} \ k=1,\ldots,m.
\end{equation}

We will construct $\ddx$ satisfying \eqref{xddot}inductively.  To 
start the induction we set
\begin{equation*}
\begin{split}
\da^{(0)}_k &=\da_k  \  \text{for} \ k=1, \ldots, m, \\ 
\dx^{(0)}   &=\sum^m_{k=1} \da^{(0)}_k \chi_{A_k}= \dx, \\
k_0  &=\max(\{k: S(\dx^{(0)},k)\neq \emptyset\}\cup \{0\}).
\end{split}
\end{equation*}

If $k_0 =0$ then $\dx^{(0)}$ satisfies \eqref{xddot} and we are 
done.  If $k_0 > 0$ then
$$
S(\dx^{(0)},k)= \emptyset \ \text{for} \ k>k_0.
$$

In the inductive process we will describe a way of defining a 
sequence of non-ne\-ga\-tive integers $k_0 > k_1 >k_2> \ldots$ 
and a sequence of simple functions 
$(\dx^{(0)},\dx^{(1)},\dx^{(2)}, \ldots)$ 
 so that for all $n$
$$
S(\dx^{(n)},k)= \emptyset \ \text{for} \ k>k_n.
$$

Once these sequences are defined we observe that since $k_0 \leq 
m$ and the sequence $(k_n)_n$ is a strictly decreasing sequence 
of non-ne\-ga\-tive integers, there exists $N \leq m+1$, so that 
$k_n =0$ and $\dx^{(N)}$ satisfies \eqref{xddot}.

To describe the inductive process, suppose that 
$(\dx^{(\nu)})^n_{\nu=0}$ and $k_0>k_1>\ldots >k_n>0$ have been 
defined so that for all $\nu \leq n$:
\begin{eqnarray}
&\ \dx^{(\nu)} = \sum^m_{k=1}\da_k^{(\nu)}\chi_{A_k},& \no \\
&\  S(\dx^{(\nu)}, k_\nu) \neq \emptyset&
\ \text{if} \ k_\nu > 0,      \lb{knu} \\
&\ S(\dx^{(\nu)},k) = \emptyset &\ \text{for} \ k>k_\nu,  \lb{bigger}  \\
&\ {\disp \frac{\da^{(\nu)}_k}{\da_k} \in \left[1, 
\frac{3}{3-\de}\right)} &\ \text{for} \ 
k=1,\ldots, m,  \lb{delta} \\
&\ \da_k^{(\nu)}=\da_k& \ \text{for} \ k \notin 
\bigcup^{\nu-1}_{\alpha=0}
 S(\dx^{(\alpha)}, k_\alpha),  \lb{union}  \\
&\ |\da_k|=|\da_l| \Longrightarrow 
|\da_k^{(\nu)}|=|\da_l^{(\nu)}|,&  
\lb{new}  \\
&\ |\da_{k_\nu}| > |\da_{k_{\nu-1}}|& \ \text{if} \  k_\nu>0.   
\lb{growth}  
\end{eqnarray}

Now define
\begin{equation*}
\da_j^{(n+1)}=
\begin{cases}
\sgn (\da_j)\cdot3\cdot|\da_{k_{n}}^{(n)}|  &\text{if $j \in S (\dx^{(n)},k_n)$,}\\
\da_j^{(n)} &\text{if $j \notin S(\dx^{(n)},k_n)$,}\\
\end{cases}
\end{equation*}
\begin{equation*}
\begin{split}
\dx^{(n+1)} &=\sum^m_{k=1}\da_k^{(n+1)} \chi_{A_k}, \\
k_{n+1} &= \max(\{k:S(\dx^{(n+1)},k)\neq \emptyset\} \cup \{0\}).
\end{split}
\end{equation*}

To prove the induction step we need to show that 
\eqref{knu}-\eqref{growth} are satisfied for $\nu=n+1$.

Clearly, if $k_{n+1} \neq 0$ then $S(\dx^{(n+1)},k_{n+1})\neq 
\emptyset$ so \eqref{knu} holds for $\nu = n+1$.  Similarly 
\eqref{bigger} holds for $\nu=n+1$ by the definition of $k_{n+1}$.

To verify \eqref{delta} for $\nu=n+1$ we first observe that if $ 
j \notin S(\dx^{(n)},k_n)$ then, by definition, 
$\da^{(n+1)}=\da^{(n)}$ and by \eqref{delta} we get
\begin{equation}\lb{deltakn}
\frac{\da_j^{(n+1)}}{\da_j} = \frac{\da_j^{(n)}}{\da_j} \in 
\left[1,\frac{3}{3-\de}\right) \ \text{for} \ j \notin S 
(\dx^{(n)}, k_n).
\end{equation}

Thus it only remains to check that \eqref{delta} is valid for 
$\nu=n+1$ and $j \in S (\dx^{(n)},k_n)$.  For this we first 
establish that:
\begin{equation}\lb{clean}
\da_j^{(n)} = \da_j \ \text{for} \ j \in S(\dx^{(n)}, k_n).
\end{equation}

To prove \eqref{clean}, by \eqref{union}, it is enough to show 
that for all $\alpha,\ 0 \leq \alpha < n$,
$$
S(\dx^{(n)}, k_n)\cap S(\dx^{(\alpha)}, k_\al) = \emptyset.
$$

Suppose, for contradiction, that there exists $\alpha,\ 0 \leq 
\alpha < n$ and $i,\ 1 \leq i \leq m$ so that: 
$$
i \in S(\dx^{(n)}, k_n)\cap S(\dx^{(\alpha)},k_\alpha).
$$

Set $t_1=t_2=|\da_i|,\ \tilde{t}_1 = |\da_i^{(n)}|, \ 
\tilde{t}_2=|\da_i^{(\alpha)}|, \ t_3=|\da_{k_n}|, \ \tilde{t}_3 
= |\da_{k_n}^{(n)}|,\  t_4 = | \da_{k_\alpha}|,\ \tilde{t}_4 = 
|\da_{k_\alpha}^{(\alpha)}|$. Then by \eqref{delta} and 
Lemma~\ref{close}$(i)$:
$$
\frac{|\da_{k_n}|}{|\da_{k_\alpha}|}=\frac{t_3}{t_4} 
\in\left(\frac{(3-\de)^3}{3^3},\frac{3^3}{(3-\de)^3}\right) = 
\left(\frac{1}{1+ \eta},1+\eta\right).
$$
Thus, by \eqref{eta}, $|\da_{k_n}|=|\da_{k_\alpha}|$ which 
contradicts \eqref{growth}. Hence \eqref{clean}is proven.

Next, by \eqref{clean} and by definition of $S (\dx^{(n)},k_n)$ 
we see that for $j \in S(\dx^{(n)}, k_n)$:
$$
\frac{|\da_j|}{|\da_{k_n}^{(n)}|}=\frac{|\da_j^{(n)}|}{|\da_{k_n}^{(n)}|} 
\in (3-\de,3).
$$
Thus, by Lemma~\ref{close}$(ii)$ with $t_1 = |\da_j|,\ t_2 = 
|\da_{k_n}^{(n)}|,\ t_3 = |\da_j^{(n+1)}|$, we obtain
\begin{equation*}
\frac{\da_j^{(n+1)}}{\da_j}=\frac{|\da_j^{(n+1)}|}{|a_j|}=\frac{t_3}{t_1}\in 
\left(1,\frac{3}{3-\de}\right) \ \ \  \text{for}\ \ \ j \in 
S(\dx^{(n)},k_n).
\end{equation*}
Together with \eqref{deltakn} this ends the proof that 
\eqref{delta} is satisfied for $\nu = n+1$.

Next we check that \eqref{union} is valid for $\nu=n+1$, i.e.
$$
\da_k^{(n+1)}=\da_k \ \ \ \  \text{for}  \ \ \  k \notin 
\bigcup^n_{\alpha=0} S(\dx^{(\alpha)},k_\al).
$$

Let $k \notin \bigcup^n_{\alpha = 0} S (\dx^{(\alpha)}, 
k_\alpha)$.  Then $k \notin S (\dx^{(n)}, k_n)$ and, by 
definition, $\da_k^{(n+1)}=\da_k^{(n)}$. And since $k \notin 
\bigcup^{n-1}_{\alpha = 0} S (\dx^{(\alpha)},
 k_\alpha)$, by \eqref{union}, $\da_k^{(n)}= \da_k$.  
 Thus \eqref{union} holds for $\nu =n+1$.

Our next step is to check \eqref{new} for $\nu = n+1$.  We know, 
by \eqref{new}, that if $|\da_k|=|\da_l|$ then $|\da_k^{(n)}| = 
|\da_l^{(n)}|$.  Thus $k \in S(\dx^{(n)}, k_n)$ if and only if $l 
\in S(\dx^{(n)}, k_n)$.  In either case it follows directly from 
the definition that $|\da_k^{(n+1)}| = |\da_l^{(n+1)}|$, i.e. 
\eqref{new} holds for $\nu=n+1$.

Our final step is to verify \eqref{growth} for $\nu=n+1$, i.e. to 
show that if $k_{n+1}>0$ then
$$
|\da_{k_{n+1}}|>|\da_{k_n}|.
$$

Since $(|\da_k|)^m_{k=1}$ are arranged in a non-increasing order 
and $k_{n+1}=\max(\{k:S(\dx^{(n+1)},k)\neq \emptyset\}\cup 
\{0\})>0$ it is enough to prove that
\begin{equation}\lb{last}
|\da_k| \leq |\da_{k_n}| \Longrightarrow S(\dx^{(n+1)}, k)= 
\emptyset .
\end{equation}

If $|\da_k|=|\da_{k_n}|$ then by \eqref{new} for $\nu=n+1$ we get 
$|\da_k^{(n+1)}| = |\da_{k_n}^{(n+1)}|$ and thus
\begin{equation}\lb{last1}
S(\dx^{(n+1)},k)=S(\dx^{(n+1)}, k_n).
\end{equation}
Notice that $k_n \notin S (\dx^n, k_n)$ so 
$\da_{k_n}^{(n+1)}=\da_{k_n}^{(n)}$. Hence, when $j \notin 
S(\dx^{(n)},k_n)$ we get $\da_j^{(n+1)}=\da_j^{(n)}$ and
$$
\frac{|\da_j^{(n+1)}|}{|\da_{k_n}^{(n+1)}|}= 
\frac{|\da_j^{(n)}|}{|\da^{(n)}_{k_n}|} \notin (3-\de,3).
$$
Thus $j \notin S(\dx^{(n+1)}, k_n)$.

If $j\in S (\dx^{(n)}, k_n)$ then, by definition,
$$
\frac{|\da_j^{(n+1)}|}{|\da_{k_n}^{(n+1)}|}=\frac{|\da_j^{(n+1)}|}{|\da^{(n)}_{k_n}|} 
=3 \notin (3-\de,3).
$$
Hence $S(\dx^{(n+1)}, k_n)= \emptyset$ and by \eqref{last1} we 
see that 
\begin{equation}\lb{last2}
|\da_k|=|\da_{k_n}| \Longrightarrow S(\dx^{(n+1)},k)= \emptyset .
\end{equation}

Next we consider the case:
$$
|\da_k| < |\da_{k_n}|.
$$
In this case $k>k_n$ and by definition of $k_n,\ S(\dx^{(n)},k)= 
\emptyset$.

By \eqref{delta}
$$
\frac{|\da_k^{(n)}|}{|\da_k|} \in 
\left[1,\frac{3}{3-\de}\right),\ \ 
\frac{|\da^{(n)}_{k_n}|}{|\da_{k_n}|} \in 
\left[1,\frac{3}{3-\de}\right).
$$

Thus if we set $t_1=|\da_k|,\ \tilde{t}_1=|\da^{(n)}_k|,\ t_2 = 
|\da_{k_n}|,\ \tilde{t}_2=|\da^{(n)}_{k_n}|$, by 
Lemma~\ref{close}$(iii)$ we obtain
$$
\frac{|\da_k^{(n)}|}{|\da^{(n)}_{k_n}|}=\frac{\tilde{t}_1}{\tilde{t}_2}<\frac{9}{(3-\de)^2} 
< 3-\de.
$$
Thus $k \notin S (\dx^{(n)},k_n)$ and, by definition, 
$\da_k^{(n+1)}=\da^{(n)}_k$. Further, for all $j \notin 
S(\dx^{(n)},k_n)$ we have $\da_j^{(n+1)}=\da_j^{(n)}$, and 
therefore
$$
\frac{|\da_j^{(n+1)}|}{|\da_k^{(n+1)}|}=\frac{|\da_j^{(n)}|}{|\da^{(n)}_k|} 
\notin (3-\de,3),
$$
since $S(\dx^{(n)},k)=\emptyset$.  Hence
$$
S(\dx^{(n+1)},k) \subset S(\dx^{(n)}, k_n).
$$
But if $j \in S(\dx^{(n+1)},k)\cap S(\dx^{(n)}, k_n)$ then
$$
\frac{|\da_j^{(n+1)}|}{|\da_{k_n}^{(n)}|}=3,\ \ \   
\frac{|\da_j^{(n+1)}|}{|\da_k^{(n)}|} 
=\frac{|\da_j^{(n+1)}|}{|\da_k^{(n+1)}|} \in (3-\de,3),
$$ 
and if we set $t_1 =|\da_j^{(n+1)}|, \ t_2 = |\da_{k_n}|, \ 
\tilde{t}_2=|\da_{k_n}^{(n)}|,\ t_3=|\da_k|,\ 
\tilde{t}_3=|\da_k^{(n)}|$ then by \eqref{delta} and by 
Lemma~\ref{close}$(iv)$ we get
$$
\frac{|\da_{k_n}|}{|\da_k|}= \frac{t_2}{t_3}\in 
\left(\frac{(3-\de)^2}{9},\frac{3}{3-\de}\right) \subset 
\left(\frac{1}{1+\eta},1+\eta\right).
$$
Thus by \eqref{eta}, $|\da_{k_n}| = |\da_k|$ which contradicts 
our assumption that $|\da_{k_n}|>|\da_k|$.

Thus $S(\dx^{(n+1)},k)=\emptyset$ if $|\da_k|<|\da_{k_n}|$, which 
together with \eqref{last2} concludes the proof of \eqref{last} 
and \eqref{growth} for $\nu=n+1$.

Note that \eqref{growth} for $\nu = n+1$, implies that $k_{n+1}< 
k_n$.

This ends the proof of the inductive process.

To finish the proof of the proposition we notice, as indicated 
above, that since $(k_n)_{n\geq 0}$ is a strictly decreasing 
sequence of nonnegative integers, it must be finite, i.e. there 
exists $N\leq m+1$, so that $k_N = 0$.

Set $\ddx = \dx^{(N)}$. Then, by \eqref{bigger},
$$
S(\dx^{(N)}, k) = \emptyset,\ \ \  \ \text{for all} \ k >k_N =0,
$$
and, by \eqref{delta},
\begin{equation*}
\begin{split}
\| \ddx - \dx \| &=\big\| \sum^m_{k=1}(\da_k^{(N)} - \da_k) 
\chi_{A_k} \big\| = \big\| \sum^m_{k=1} \da_k 
\left(\frac{\da^{(N)}_k}{\da_k}-1\right) 
\chi_{A_k} \big\| \\
&\leq \left(\frac{3}{3- \de}-1\right) \cdot
 \big\| \sum^m_{k=1} \da_k \chi_{A_k} \big\|
 = \frac{\de}{3-\de} \| \dx \|.
\end{split}
\end{equation*}

Thus, by \eqref{xdot}, when $\ \de < 1/8$ we have:
\begin{equation*}
\begin{split}
\| \ddx - x \| &\leq \| \ddx - \dx \| + \| \dx - x \| \leq 
\frac{\de}{3-\de} \| \dx \| + \eta \| x \| \leq 
\frac{\de}{3-\de}\cdot (1+\eta) \|x\| + \eta \|x\| \\& = 
\left[\frac{\de}{3-\de}\cdot\left(\frac{3}{3-\de}\right)^3 + 
\left(\frac{3}{3-\de}\right)^3-1\right] \| x \| = 
\frac{81-(3-\de)^4}{(3-\de)^4}\ \|x \|\\& < \frac{3}{2}\ \de 
\|x\|.
\end{split}
\end{equation*}

Finally note that, by \eqref{delta}, {for all}  $j=1,\dots,m$,
$$ \frac{\da_j^{(N)}}{\da_j}\ge 1. $$ 

Combining the above two inequalities and   Lemma~\ref{step1} we 
get
$$\left(1+\frac32 \de\right)\|x\|>\|\ddx\|=\|\dx^{(N)}\|\ge\|\dx\|\ge\|x\|.$$

This ends the proof of the proposition.
\end{proof}

\section{Remarks about duality and narrow operators}\lb{dual}

As mentioned above (see Remark~\ref{p<2}) our proof does not work in the
case $1\le p<2$. Unfortunately this case cannot be handled by 
duality arguments either. Indeed 
the dual of $L_{w,p}$ is never
isometric to a Lorentz space, even if we restrict ourselves to the more
classical spaces $L_{p,q}$ (cf. e.g. \cite{BSh}).
The isometric structure of the  duals of  Lorentz spaces $L_{w,p}$ 
has  been described in \cite{H53}, however this description is so 
complicated that there are later papers describing the isomorphic
structure of   duals of  Lorentz spaces $L_{w,p}$, see 
\cite{A78,R82,S90,KMa,KM} and their references. These results require too
much technical notation to be stated here, we just mention the result for 
spaces $L_{p,q}$.

Recall that the space $L_{p,q}$ for $p,q\in(1,\infty)$ is defined as
$L_{q,w}$, where $w(t)= qt^{p/q-1}/p$. It is easily seen that $\|\cdot\|_{q,w}$
with this weight is not a norm when $q>p$, since then $w(t)$ is decreasing
rather than decreasing. Nevertheless $L_{p,q}$ is a linear space also
for $q>p$ and it can be shown that $L_{p,q}$ can be made into a Banach space
when $p>1$ by introducing an actual norm $|\!|\!|\cdot|\!|\!|_{p,q}$
which satisfies $\|f\|_{p,q}\le |\!|\!|f|\!|\!|_{p,q}\le C(p,q)\|f\|_{p,q}$.

For spaces $L_{p,q}$ we have the following isomorphic characterization
of the dual space
$$(L_{p,q})^*\cong L_{p',q'},$$
where $1/p+1/p'=1/q+1/q'=1$.

Another difficulty with the duality approach comes from
the fact  that in general the 
conjugate operator $T^*$ to a narrow operator $T:E\lra E$ 
need not be narrow (for any r.i. space $E$ with $E^*$ having an 
absolutely continuous norm to consider the notion of narrow 
operators) \cite[p. 60]{PP90}. 
In general we can only conclude that if $P$ is
 a projection onto a finite-codimensional subspace $X 
\subseteq E$ then $P^*$ is also a projection onto a 
finite-codimensional subspace in $E^*$. And since $\|P\| = 
\|P^*\|$, we have the following obvious

\begin{prop}\lb{1} Suppose that for a reflexive r.i. space $E$ there 
is a constant 
$k_E$ such that $\|P\| \geq k_E$ for any finite-codimensional 
projection $P$ in $E$. Then the same is true for 
finite-codimensional projections in $E^*$.
\end{prop}

Thus we obtain the immediate:
\begin{cor}
\begin{itemize}
\item[(1)] 
For any $p>2$ there exists a constant $k_p>1$ so that for 
every projection $P$ from $(L_{w,p})^*$ onto a finite-codimensional 
subspace we have $\|P\| \geq k_p$.
\item[(2)] For any $q$, $1< q<2$ and any $p\in(1,\infty)$
there exists an
equivalent norm   $|\!|\!|\cdot|\!|\!|$ on $L_{p,q}$,   so that  
there exists a constant $k_q>1$ such that for 
every projection $P$ from $L_{p,q}$ onto a finite-codimensional 
subspace we have $|\!|\!|P|\!|\!| \geq k_q$.
\end{itemize}
\end{cor}

We recall here that it follows from \cite[Theorem~4]{pams} 
(Theorem~\ref{pams}
above) that in $L_{w,p}$ with the usual norm, for all $p$, $1\le p<\infty$ 
and all weights $w$, we have $\|P\| >1$ for all finite-codimensional
projections.

In the remainder of this section we study operators which 
are conjugate to narrow. 
We start with the definition:
\begin{defn} \lb{sn} We say that an operator $T \in {\cal L}(E)$ on an r.i. 
space $E$ with absolutely continuous norm is {\it *-narrow} 
provided $T^* \in {\cal L}(E^*)$ is narrow.
\end{defn}

Next we will need the 
following notion.

\begin{defn} Let $E$ be an r.i. space with absolutely continuous norm. A 
subset $M \subseteq E$ we call {\it poor} if for each measurable 
subset $A \subseteq [0,1]$ and each $\varepsilon >0$ there exists 
a decomposition $A = A^+ \bigsqcup A^-$ into measurable subsets 
such that
$$
| \int\limits_{A^+} f d\mu - \int\limits_{A^-} f d\mu | < 
\varepsilon
$$
for every $f \in M$.
\end{defn}

Denote by $B(X)$ the closed unit ball of a Banach space $X$.
Here and in the sequel by ${\cal L}(X,Y)$ we denote 
the space of all continuous linear operators acting from $X$ to 
$Y$ and use ${\cal L}(X)$ instead of ${\cal L}(X,X)$. 

\begin{prop}\lb{2} Let $E$ be a reflexive r.i. space. An operator 
$T \in {\cal 
L}(E)$ is narrow if and only if $T^*(B(E^*))$ is a poor subset of 
$E^*$.
\end{prop}

\begin{proof} For each $x \in E$ we have
$$
\|Tx\| = \sup\limits_{f \in B(E^*)} |f(Tx)| = \sup\limits_{f \in 
B(E^*)} |(T^*f)(x)|.
$$

Now suppose that $x^2 = \chi_{_A}$ and put $A^+ = x^{-1}(1)$ and 
$A^- = x^{-1}(-1)$. Then
$$
\|Tx\| = \sup\limits_{f \in B(E^*)} | \int\limits_{A^+} T^*f d\mu 
- \int\limits_{A^-} T^*f d\mu |.
$$

This observation completes the proof.
\end{proof}

As an immediate corollary we obtain:

\begin{cor} Let $E$ be a reflexive r.i. space with absolutely 
continuous norm. An operator 
$T \in {\cal 
L}(E)$ is *-narrow if and only if $T(B(E))$ is a poor subset of 
$E$.
\end{cor}

\begin{proof}
Indeed, $T \in {\cal L}(E)$ is *-narrow if and only if $T^* \in 
{\cal L}(E^*)$ is narrow, and 
since $T^{**} = T$ for operators $T$ on a reflexive space,
by Proposition~\ref{2},  we have $T^* \in {\cal L}(E^*)$ is *-narrow
if and only if $T^{**}(B(E^{**}))=T(B(E))$ is poor.
\end{proof}

Denote by $Narr(E)$ and $*-Narr(E)$ the sets of all narrow and 
respectively *-narrow operators on $E$.

The existence of a narrow operator whose conjugate is not narrow 
shows that the both inclusions $Narr(E) \subseteq *-Narr(E)$ and 
$Narr(E) \supseteq *-Narr(E)$ are false. Indeed, let $T$ be 
narrow on $E$ and $T^*$ not be narrow on $E^*$. If $T$ were 
*-narrow, then $T^*$ would be narrow. Thus, first inclusion 
fails. Let $T^*$ be narrow in $E^*$ and $T$ not be. Then $T$ is 
*-narrow in $E$ and the second inclusion fails.

The intersection $Narr(E) \bigcap *-Narr(E)$ contains the set 
${\cal K}(E)$ of all compact operators on $E$. The inclusion 
${\cal K}(E) \subseteq Narr(E)$ is very simple and proved in 
\cite[p. 55]{PP90}. The inverse inclusion is also simple and 
follows from

\begin{prop}\lb{3} Let $E$ be an r.i. space with absolutely 
continuous norm. 
Then every relatively compact set $K \in E$ is poor.
\end{prop}
\begin{proof} Given $\varepsilon >0$ and a measurable subset $A 
\subseteq [0,1]$. Pick any $\frac{\varepsilon}{2}$-net in $K$, 
say $\{f_1,...,f_m\}$. Let $\{r_n(A)\}$ be any ``Rademacher'' 
type sequence on $A$, i.e. $(r_n(A))^2 = \chi_{_A}$ and $r_n(A) 
\stackrel{w}{\longrightarrow} 0$. Now choose $n$ so that $|\int 
f_i r_n(A) d\mu|< \frac{\varepsilon}{2}$ for each $i=1,...,m$. 
Then for each $f \in K$ we obtain
$$
|\int f r_n(A) d\mu | \leq |\int f_i r_n(A) d\mu | + \|f_i - f\| 
< \varepsilon
$$
for suitable $i$. The proposition is proved.
\end{proof}

A principle difference between the notions of narrow and *-narrow 
operators is contained in the following: *-narrowness of an 
operator is a property of the image under it of the unit ball. 
But the notion of narrow operators cannot be formulated in terms 
of the image. This fact is a consequence of the following one.

\begin{theorem}\lb{4} Let $E$ be an r.i. space with an absolutely continuous 
norm. Then there exists a complemented subspace $E_0$ of $E$ 
isomorphic to $E$ and for which there are two projections onto 
$E_0$, one of which is narrow and second is not narrow at any 
measurable subset $B \subseteq [0,1]$ of positive measure.
\end{theorem}

(An operator $T$ on an r.i. space $E$ is said to be narrow at a 
measurable subset $B \subseteq [0,1]$ if the operator $T(\bullet 
\cdot \chi_{_B})$ is narrow.)

In other words, there are two decompositions $E = E_0 \oplus E_1 
= E_0 \oplus E_2$ with rich $E_1$ and $E_2$ which is not rich at 
any measurable subset $B \in [0,1]$.

\begin{proof} First construct two subspaces $E_0$ and $E_1$. For 
convenience, we consider $E$ at $[0,1] \times [0,1]$ instead of 
$[0,1]$ (cf. \cite[Example 1, p. 57]{PP90}). We consider elements 
of $E$ as equivalence classes of functions of two variables 
$x(t_1,t_2)$. Put
\begin{equation*}
\begin{split}
E_0 &= \{  x \in E: x(t_1,t_2)= 
\int\limits_{[0,1]} x(t_1,s) ds: x\in E \; \; (a.e.) \},\\ 
 E_1 &= \{ x \in E: \int\limits_{[0,1]} x(t_1,t_2) dt_2 = 0 
\; \; (a.e.) \}.
\end{split}
\end{equation*}

So, $E = E_0 \oplus E_1$. In \cite[p. 57]{PP90} it is proved that 
$E_1$ is rich and therefore the projection $P$ of $E$ onto $E_0$ 
with $\ker P = E_1$ is narrow.

Now construct another complement $E_2$ to $E_0$. We use the 
following simple argument. Suppose $Z = X \oplus Y$ where $X,Y,Z$ 
are Banach spaces and $T \in {\cal L}(Y,X)$. Then $Z = X \oplus 
\{y+Ty:y \in Y \}$ is a decomposition into closed subspaces (it 
is easily verified). For $T \neq 0$ we somewhat ``bend'' the 
subspace $Y$. Our purpose is to ``bend'' the subspace $E_1$ by a 
suitable $T$.

Put $D_0 = [0, \frac{1}{2}) \times [0,1], \; \; \; D_1 = 
[\frac{1}{2},1] \times [0,1], \; \; \; E_i' = \chi_{_{D_0}}E_i, 
\; \; \; E_i'' = \chi_{_{D_1}}E_i$ for $i=0,1$. So we have $E_i = 
E_i' \oplus E_i'', \; \; \; i=0,1$. Let $T$ be an isomorphic 
embedding of $E_1$ into $E_0$ with $TE_1' \subseteq E_0''$ and 
$TE_1'' \subseteq E_0'$ (such an operator exists since both 
$E_0'$ and $E_0''$ are evidently isomorphic to $E$).

Next we show that the subspace $E_2 = \{y+Ty: y \in E_1 \}$ (which is a 
complement to $E_0$) is not rich at any subset $B \subseteq [0,1] 
\times [0,1]$. Fix any such $B$ of positive measure. Then at 
least one of the two subsets $B \bigcap D_1$ or $B \bigcap D_2$ 
is of positive measure, say first of them. Denote it by $A$. Let 
$x$ be any element of $E$ with $x^2 = \chi_{_A}$ and $y+Ty$ be 
any element of $E_2$ ($y \in E_1$). Since $y \chi_{_{D_1}}$ and 
$y \chi_{_{D_2}}$ belong to $E_1$ we may write:
$$
\|x-y-Ty\|=\|x-y\chi_{_{D_1}}-y\chi_{_{D_2}}-T(y\chi_{_{D_1}})-T(y\chi_{_{D_2}})\|\geq
$$
(since the restriction of a function to a subset of the domain is 
a projection of norm one)
$$
\max \{\|x-y\chi_{_{D_1}}-T(y\chi_{_{D_2}})\|, \; 
\|y\chi_{_{D_2}}+T(y\chi_{_{D_1}})\| \} =
$$
$$
\max \{\|\int x dt_2 - T(y\chi_{_{D_2}})+(x- \int x 
dt_2)-y\chi_{_{D_1}}\|, \; \|y\chi_{_{D_2}}+T(y\chi_{_{D_1}})\| 
\} \geq
$$
$$
\max \{\| \int x dt_2 - T(y\chi_{_{D_2}})\|, \; \alpha^{-1} \|x - 
\int x dt_2 - y\chi_{_{D_1}} \|, \; \alpha^{-1} \|y\chi_{_{D_2}} 
\|, \; \|T(y\chi_{_{D_1}}) \| \},
$$
where $\alpha$ is the norm of the projection of $E$ onto $E_1$ 
having the kernel $E_0$ (note that the projection of $E$ onto 
$E_0$ with kernel $E_1$ is a   conditional expectation operator
and therefore is of norm one).

Now suppose to the contrary that $E_2$ is rich at $B$. Then for 
each $\varepsilon > 0$ there are $x_{\varepsilon} \in E$ with 
$x_{\varepsilon}^2 = \chi_{_A}$ and $y_{\varepsilon} \in E_1 \; 
\; \; (y_{\varepsilon}+Ty_{\varepsilon} \in E_2)$ such that 
$\|x_{\varepsilon}-y_{\varepsilon}-Ty_{\varepsilon} \| < 
\varepsilon$.

By the above estimates we have
$$
\|\int x_{\varepsilon} dt_2 - T(y_{\varepsilon}\chi_{_{D_2}})\| < 
\varepsilon, \ \ \ \|x_{\varepsilon}-\int x_{\varepsilon} 
dt_2 - y_{\varepsilon}\chi_{_{D_1}} \| < \alpha \varepsilon, 
 \ \ \   \|y_{\varepsilon} \chi_{_{D_2}} \| < \alpha 
\varepsilon,  \ \ \  \|T(y_{\varepsilon}\chi_{_{D_1}})\| < 
\varepsilon
$$
whence
$$
\|\chi_{_A}\|=\|x_{\varepsilon}\|=\|\int x_{\varepsilon} dt_2 + 
(x_{\varepsilon} - \int x_{\varepsilon} dt_2)\| \leq
$$
$$
\| \int x_{\varepsilon} dt_2 - T(y_{\varepsilon} \chi_{_{D_2}})\| 
+ \|T\| \|y_{\varepsilon}\chi_{_{D_2}}\| + \|x_{\varepsilon} - 
\int x_{\varepsilon} dt_2 - y_{\varepsilon}\chi_{_{D_1}} \| + 
\|y_{\varepsilon} \chi_{_{D_1}} \| \leq
$$
$$
\varepsilon (1 + \alpha \|T\| + \alpha + \|T^{-1}\|).
$$

The contradiction completes the proof.
\end{proof}

Theorem~\ref{4} implies a negative answer to the following 
question of A. Plichko and M. Popov \cite[Question 2, p. 71]{PP90}:

Given an r.i. space $E, \;\;\; E \neq L_2$ such that $E^*$ has 
absolutely continuous norm. Suppose that $E = X \oplus Y$ and 
$E^* = Y^{\bot} \oplus X^{\bot}$ where $X^{\bot}\ ( Y^{\bot})$ 
consists of all functionals vanishing on $X$ (on $Y$). Is $X \in 
E$ rich if and only if $Y^{\bot} \in E^*$ is?

Indeed, we have constructed decompositions $E = E_0 \oplus E_1 = 
E_0 \oplus E_2$ with $E_1$ rich and $E_2$ not rich, 
which imply the dual decompositions $E^* = 
E_1^{\bot} \oplus E_0^{\bot} = E_2^{\bot} \oplus E_0^{\bot}$. If 
we assume an affirmative answer to the above question then 
$E_0^{\bot}$ should be rich and therefore $E_2$ also should be 
rich.


\begin{thebibliography}{10}

\bibitem{A78}
{\sc G.~D. Allen}, {\em Duals of {L}orentz spaces}, Pacific J. 
Math., 77
  (1978), pp.~287--291.

\bibitem{ABJS79}
{\sc S.~Axler, I.~D. Berg, N.~Jewell, and A.~Shields}, {\em 
Approximation by
  compact operators and the space ${H}\sp{\infty }+{C}$}, Ann. of Math. (2),
  109 (1979), pp.~601--612.

\bibitem{BO87}
{\sc H.~Bang and E.~Odell}, {\em On the best compact 
approximation problem for
  operators between ${L}\sb p$-spaces}, J. Approx. Theory, 51 (1987),
  pp.~274--287.

\bibitem{BSh}
{\sc C.~Bennett and R.~Sharpley}, {\em Interpolation of 
operators}, Academic
  Press Inc., Boston, MA, 1988.

\bibitem{BL85}
{\sc Y.~Benyamini and P.~K. Lin}, {\em An operator on ${L}\sp p$ 
without best
  compact approximation}, Israel J. Math., 51 (1985), pp.~298--304.

\bibitem{F90}
{\sc C.~Franchetti}, {\em The norm of the minimal projection onto 
hyperplanes
  in ${L}\sp p[0,1]$ and the radial constant}, Boll. Un. Mat. Ital. B (7), 4
  (1990), pp.~803--821.

\bibitem{F92}
\leavevmode\vrule height 2pt depth -1.6pt width 23pt, {\em Lower 
bounds for the
  norms of projections with small kernels}, Bull. Austral. Math. Soc., 45
  (1992), pp.~507--511.

\bibitem{FS95}
{\sc C.~Franchetti and E.~M. Sem{\"e}nov}, {\em A {H}ilbert space
  characterization among function spaces}, Anal. Math., 21 (1995), pp.~85--93.

\bibitem{H53}
{\sc I.~Halperin}, {\em Function spaces}, Canadian J. Math., 5 
(1953),
  pp.~273--288.

\bibitem{JMST}
{\sc W.~B. Johnson, B.~Maurey, G.~Schechtman, and L.~Tzafriri}, 
{\em Symmetric
  structures in {B}anach spaces}, Mem. Amer. Math. Soc., 19 (1979), pp.~v+298.

\bibitem{KP96}
{\sc V.~M. Kadets and M.~M. Popov}, {\em The {D}augavet property 
for narrow
  operators in rich subspaces of the spaces ${C}[0,1]$ and ${L}\sb 1[0,1]$},
  Algebra i Analiz, 8 (1996), pp.~43--62.
\newblock (Russian), translation in {\it St. Petersburg Math. J.} {\bf 8}
  (1997), 571--584.

\bibitem{KSSW00}
{\sc V.~M. Kadets, R.~V. Shvidkoy, G.~G. Sirotkin, and 
D.~Werner}, {\em Banach
  spaces with the {D}augavet property}, Trans. Amer. Math. Soc., 352 (2000),
  pp.~855--873.

\bibitem{KSW}
{\sc V.~M. Kadets, R.~V. Shvidkoy, and D.~Werner}, {\em Narrow 
operators and
  rich subspaces of {B}anach spaces with the {D}augavet property}.
\newblock preprint.

\bibitem{KR}
{\sc N.~J. Kalton and B.~Randrianantoanina}, {\em Surjective 
isometries of
  rearrangement-invariant spaces}, Quart. J. Math. Oxford, 45 (1994),
  pp.~301--327.

\bibitem{KMa}
{\sc A.~Kami\'nska and L.~Maligranda}, {\em On {L}orentz spaces
  {$\Gamma_{p,w}$}}.
\newblock preprint.

\bibitem{KM}
{\sc A.~Kami\'nska and M.~Masty{\l}o}, {\em Duality and classical 
operators in
  function spaces}.
\newblock preprint.

\bibitem{LT2}
{\sc J.~Lindenstrauss and L.~Tzafriri}, {\em Classical {B}anach 
spaces, Vol. 2,
  {F}unction spaces}, Springer--Verlag, Berlin--Heidelberg--New York, 1979.

\bibitem{Lor50}
{\sc G.~G. Lorentz}, {\em {Some new functional spaces}}, Ann. of 
Math., 51
  (1950), pp.~37--55.

\bibitem{Lor51}
\leavevmode\vrule height 2pt depth -1.6pt width 23pt, {\em On the 
theory of
  spaces {$\Lambda$}}, Pacific J. Math., 1 (1951), pp.~411--429.

\bibitem{L66}
{\sc G.~Y. Lozanovski{\u\i}}, {\em On almost integral operators in
  ${K}{B}$-spaces}, Vestnik Leningrad. Univ., 21 (1966), pp.~35--44.
\newblock (Russian).

\bibitem{O}
{\sc T.~Oikhberg}, {\em The {D}augavet property of 
{$C^*$}-algebras and
  non-commutative {$L_p$}-spaces}, Positivity.
\newblock to appear.

\bibitem{PP90}
{\sc A.~M. Plichko and M.~M. Popov}, {\em Symmetric function 
spaces on atomless
  probability spaces}, Dissertationes Math. (Rozprawy Mat.), 306 (1990),
  pp.~1--85.

\bibitem{Popov87}
{\sc M.~M. Popov}, {\em On norms of projectors in ${L}\sb p(\mu)$ 
with
  ``small'' kernels}, Funktsional. Anal. i Prilozhen., 21 (1987), pp.~84--85.

\bibitem{pams}
{\sc B.~Randrianantoanina}, {\em Contractive projections in 
nonatomic function
  spaces}, Proc. Amer. Math. Soc., 123 (1995), pp.~1747--1750.

\bibitem{R82}
{\sc S.~Reisner}, {\em On the duals of {L}orentz function and 
sequence spaces},
  Indiana Univ. Math. J., 31 (1982), pp.~65--72.

\bibitem{S90}
{\sc E.~Sawyer}, {\em Boundedness of classical operators on 
classical {L}orentz
  spaces}, Studia Math., 96 (1990), pp.~145--158.

\end{thebibliography}

\def\polhk#1{\setbox0=\hbox{#1}{\ooalign{\hidewidth
  \lower1.5ex\hbox{`}\hidewidth\crcr\unhbox0}}} \def\cprime{$'$}

\end{document}